\documentclass{amsart}

\newtheorem{theorem}{Theorem}[section]
\newtheorem{lemma}[theorem]{Lemma}
\newtheorem{corollary}[theorem]{Corollary}
\newtheorem{proposition}[theorem]{Proposition}

\theoremstyle{definition}

\theoremstyle{remark}

\numberwithin{equation}{section}

\DeclareMathOperator{\re}{Re}
\DeclareMathOperator{\im}{Im}
\DeclareMathOperator{\supp}{supp}
\DeclareMathOperator{\Spec}{Spec}

\newcommand{\bbC}{{\mathbb C}}
\newcommand{\bbZ}{{\mathbb Z}}
\newcommand{\bbR}{{\mathbb R}}
\newcommand{\bbH}{{\mathbb H}}
\newcommand{\bbN}{{\mathbb N}}
\newcommand{\bbB}{{\mathbb B}}
\newcommand{\calK}{{\mathcal K}}
\newcommand{\calV}{{\mathcal V}}
\newcommand{\calA}{{\mathcal A}}
\newcommand{\calM}{{\mathcal M}}
\newcommand{\calD}{{\mathcal D}}

\newcommand{\calS}{{\mathcal S}}
\newcommand{\calE}{{\mathcal E}}
\newcommand{\calF}{{\mathcal F}}

\newcommand{\calI}{{\mathcal I}}
\newcommand{\cmet}{d}
\newcommand{\del}{\partial}
\newcommand{\cinf}{C^\infty}
\newcommand{\cinfd}{\dot C^\infty}
\newcommand{\diag}{\Lambda}
\newcommand{\dlift}{\tilde\Lambda}
\newcommand{\bdiag}{\del\Lambda}
\newcommand{\krnl}[1]{\calK_{#1}}
\newcommand{\tkrnl}[1]{\tilde\calK_{#1}}
\newcommand{\metr}{\calM_X}
\newcommand{\xx}{X\times X}
\newcommand{\sxx}{{X\mathbin{\tilde\times}X}}
\newcommand{\xU}{\mathbin{\overline\cup}}
\newcommand{\bX}{{\del X}}
\newcommand{\bxbx}{\bX\times\bX}

\newcommand{\sxxx}{\widetilde{X^3}}
\newcommand{\sxbx}{{X\mathbin{\tilde\times}\bX}}
\newcommand{\norm}[1]{\Vert #1\Vert}

\newcommand{\Sz}{{\check S_\zeta}}
\newcommand{\deriv}[1]{\frac{\del}{\del#1}}
\newcommand{\diff}{\text{Diff}}
\newcommand{\Emin}{{\mathfrak e}}

\newcommand{\hfd}{\Omega^{1/2}}
\newcommand{\phgd}[1]{{\calA}_{#1}}
\newcommand{\phgk}[1]{{\calA}_{#1}(\sxx;\> \tau\hfd)}
\newcommand{\op}[2]{\tilde\Psi^{#1}_{#2}(X;\> \hfd)}
\newcommand{\psbop}[1]{\Psi^{#1}(\bX)}
\newcommand{\vrop}[1]{\Psi_{#1}(X;\> \hfd)}
\newcommand{\cnrml}[1]{I^{#1}(\sxx, \dlift;\> \tau\hfd)}

\newcommand{\hfdhil}{L^2(X; \rho^{-(n+1)/2} \hfd)}

\begin{document}

\title[Scattering theory and deformations]{Scattering theory and 
deformations of asymptotically hyperbolic metrics}
\author[D. Borthwick]{David Borthwick}
\address{Department of Mathematics and Computer Science, 
Emory University, Atlanta}
\email{davidb@mathcs.emory.edu}
\thanks{Supported in part by NSF grant DMS-9401807 and by an NSF
Postdoctoral Fellowship.}
\date{November 20, 1997}
\begin{abstract}
For an asymptotically hyperbolic metric on the interior of
a compact manifold with boundary, 
we prove that the resolvent and scattering operators are continuous functions 
of the metric in the appropriate topologies.
\end{abstract}
\maketitle
\tableofcontents

\section{Introduction}
Let $X$ be a compact manifold with boundary and $g$ an asymptotically 
hyperbolic metric on the interior of $X$, in the sense that all sectional 
curvatures approach $-1$ at $\bX$.  We will assume that $g$ is of the form 
$$
g = \rho^{-2} h,
$$
where $\rho$  is a boundary defining function on $X$ and $h$ is a Riemannian
metric on $X$.   Such a metric $g$ is necessarily complete.  The asymptotic
curvature condition reduces to $|d\rho|_h = 1$ on $\bX$.  Any compact manifold
with boundary possesses such a metric.  The chief examples are
purely hyperbolic, with the interior of $X$ isometric to
$\bbH^{n+1}$ or a convex co-compact quotient $\bbH^{n+1}/\Gamma$.

Let $\Delta_g$ be the Laplacian associated to $g$, acting on functions.
It was proven in \cite{LP82}, \cite{Ma88} that the spectrum of 
$\Delta_g$ consists of absolutely continuous
spectrum $[n^2/4, \infty)$, plus finite point spectrum 
$\Spec(\Delta_g) \subset (0,n^2/4)$ (no embedded eigenvalues, in particular).

In \cite{MM} Mazzeo and Melrose demonstrated the
meromorphic continuation of the (modified) resolvent,
$$
R_\zeta := [\Delta_g - \zeta(n-\zeta)]^{-1},\qquad \dim X=n+1.
$$
The proof involves the careful construction of a parametrix for 
$\Delta_g - \zeta(n-\zeta)$.
This construction is part of the general program for dealing with degenerate 
elliptic boundary problems originated by Melrose and to be presented in detail
in \cite{MeBk}.  

The continuous spectrum of $\Delta_g$ may be characterized by the behavior of
generalized eigenfunctions at infinity.  
If $\re\zeta = n/2$, $\zeta \ne n/2$,
then for each $f\in\cinf(\bX)$ there exists a unique solution of 
$\Delta_g u = \zeta(n-\zeta) u$ with asymptotic behavior
$$
u = \rho^{n-\zeta} f + \rho^{\zeta}f' + O(\rho^{n/2+\epsilon}),
$$
where $f'\in\cinf(\bX)$, $\epsilon>0$.  
In fact such a solution will have a complete polyhomogeneous
expansion in $\rho$ \cite{Ma91}.

The ``scattering operator,'' defined by 
\begin{equation}\label{Sdef}
S_\zeta:f\mapsto f',
\end{equation}
is a zeroth-order pseudodifferential operator on $\bX$.  
This follows fairly directly from the results of
\cite{MM}, although the scattering operator was not considered 
there. \footnote{See \cite{MeGS} for a review of scattering theory on
manifolds with various types of regular structure at infinity.} 
For the case $X = \bbH^{n+1}/\Gamma$, the meromorphic continuation of 
the resolvent was also proven by Perry \cite{P89} using a method 
involving the scattering operator,
which was shown to be pseudodifferential explicitly.

As defined by (\ref{Sdef}), $S_\zeta$ depends on the choice of boundary 
defining function $\rho$.   To remove this dependency, we may introduce the
bundles $\Omega^{\alpha}$ of $\alpha$-densities on $\bX$, where 
$\alpha\in\bbC$.  Let $\gamma_h$ be the density on $\bX$ coming from the 
metric induced by $h$.  Then we define the normalized scattering operator
$$
\Sz : \cinf(\bX; \Omega^{1-\zeta/n}) \to \cinf(\bX, \Omega^{\zeta/n}),
$$
by
$$
\Sz: f \cdot(\gamma_h)^{1-\zeta/n} \mapsto (S_\zeta f) 
\cdot(\gamma_h)^{\zeta/n}.
$$
It is easily seen that $\Sz$ depends only on $g$ and not on $\rho$.  
For $X = \bbH^{n+1}/\Gamma$, $\cinf(\bX; \Omega^\alpha)$ 
may be conveniently realized as a space of automorphic forms on the
regular set of $\Gamma$, and the kernel of $\Sz$ may be written as an average 
over $\Gamma$ of the scattering kernel for $\bbH^{n+1}$.
This coincides with the scattering operator as
studied in \cite{Ag}, \cite{BMT}, \cite{F}, \cite{GZ97}, \cite{LP84}, 
\cite{Pa}, \cite{P89}, \cite{P95}.  

We will find it more convenient to deal with $S_\zeta$ rather than $\Sz$.
Of course the results apply to either definition.
The scattering operator also extends to a meromorphic family in $\zeta$,
and it may be derived from the resolvent by taking a certain limit at the boundary.  

The question we will address is the continuity of the resolvent
and scattering operator under deformations of the metric.  Denote by $\metr$
the space of asymptotically hyperbolic metrics on $X$, with the topology 
inherited from $\rho^{-2} \cinf(X,  T^*X\otimes T^*X)$.  This is just the $\cinf$
topology on the metric $h$.

A refinement of the construction of \cite{MM}, given in \cite{Ma91},
shows that $\tkrnl{R_\zeta}$, the lift of the Schwarz 
kernel of $R_\zeta$ to the stretched product $\sxx$, is a distribution 
with polyhomogeneous conormal singularities at the boundary.
In the notation to be introduced in \S\ref{spdsec}, 
$$
R_\zeta \in \op{-2}{\zeta,\zeta,0} + \vrop{\zeta,\zeta}.
$$
Roughly this means that the lifted
kernel has full asymptotic expansions near all boundaries, which are polyhomogeneous 
in the appropriate boundary defining function.  The subscripts give the
leading orders of these expansions.  The topology of polyhomogeneous conormal
distributions controls the behavior of all coefficients in the boundary expansions.

\begin{theorem}\label{mainthm}
For $\zeta\ne n/2$, the map
$$
\metr \ni g \mapsto R_\zeta \in \op{-2}{\zeta,\zeta,0} + \vrop{\zeta,\zeta}
$$
is continuous except at poles.
\end{theorem}

The implications for the scattering operator are easier to describe.  
From the structure of the resolvent we may deduce the local form
of the kernel:
$$
K_{S_\zeta}(y,y') = r^{-2\zeta} F(r,\theta,y) + G(y,y'),
$$
where $(y,y')$ are local coordinates for $\bxbx$,
$r=|y-y'|$, $\theta = (y-y')/r$, and $F$ and $G$ are smooth
in their respective variables.  This implies that $S_\zeta$
is a pseudodifferential operator of order $2\zeta-n$
with a one-step polyhomogeneous symbol expansion.  
From Theorem \ref{mainthm} 
we will deduce that the maps $g\mapsto F,G$ are continuous
in a $\cinf$ topology.  If $\psbop{a}$ is the space of all 
one-step polyhomogeneous pseudodifferential operators of order $a$, 
with the appropriate topology (defined in \S\ref{rsosec}), then we have the following.
\begin{theorem}\label{scthm}
For $\zeta\ne n/2$, the map
$$
\metr\ni g \mapsto S_\zeta \in \psbop{2\zeta-n}
$$ 
is continuous except at poles.
\end{theorem}

In joint work with Peter Perry, these results will be applied to study the
behavior of scattering poles and resonances (poles of the resolvent) under 
metric deformations \cite{BP}.

A converse to this result was proven in \cite{BMT}, in the special case where
$X = \bbH^3/\Gamma$.  There it was shown that
the size of a quasiconformal deformation of $\Gamma$ 
is controlled by the change in scattering operator, in the operator
topology.  The standard topology of quasiconformal deformations is, however,
essentially a $C^0$ topology, i.e. much weaker than that of $\metr$. 
In \cite{BPT}, a homeomorphism will be established between the quasiconformal 
deformation space of $\bbH^3/\Gamma$ and the space of scattering operators 
endowed with a suitably weak topology.

\bigskip\noindent{\bf Acknowledgments.}
The idea for this project arose from joint work with Ed Taylor and 
Peter Perry, and I'm indebted to them for this inspiration, as well 
as for many corrections and additions to the manuscript.
I also thank Rafe Mazzeo for some very helpful conversations, and Richard Melrose 
for some tips on asymptotic summation.  

\section{The hyperbolic model}\label{hypmod}

One of the key features of the analysis of \cite{MM} is the
relation of the general case back the model case of the ball 
$\bbB^{n+1}$ with $g$ the standard hyperbolic metric.  It is
generally more convenient to deal with the half-space model $\bbH^{n+1}$ with
coordinates $z = (x,y)$, $x\ge 0$, and the metric
$g=(dx^2+dy^2)/x^2$.

The Laplacian in these coordinates is
$$
\Delta = -(x\del_x)^2 + nx\del_x - \sum_{i=1}^n (x\del_{y_i})^2.
$$
Let $G_\zeta(z,z')$ be Schwartz kernel of the resolvent $(\Delta -
\zeta(n-\zeta))^{-1}$, with respect to the measure $dg$.
$G_\zeta$ is purely a function of the hyperbolic distance $d(z,z')$, 
given by
$$
\cosh d(z,z') = 1 + \frac{|z-z'|^2}{2xx'}.
$$
Explicitly,
$$
G_\zeta(z,z') = c_\zeta \>\Bigl(\cosh \frac{d}{2}\Bigr)^{-2\zeta} \>
F(\zeta,\zeta-\tfrac{n-1}2,2\zeta-n+1; (\cosh\tfrac{d}2)^{-2}),
$$
where $F$ is the hypergeometric function
$$
F(a,b,c;u) = 1 + \frac{a\cdot b}{1\cdot c}\>u + 
\frac{a(a+1)\cdot b(b+1)}{1\cdot 2\cdot c(c+1)}\>u^2 + \dots,
$$
and 
\begin{equation}\label{czeta}
c_\zeta = \pi^{-n/2} 2^{-2\zeta-1} 
\frac{\Gamma(\zeta)}{\Gamma(\zeta - \tfrac{n}2  + 1)}.
\end{equation}

We can define a generalized eigenfunction for $\Delta$ by taking
$$
E_\zeta(z,y') = \lim_{x'\to 0} (x')^{-\zeta} G_\zeta(z,z').
$$
$E_\zeta$ is a smooth function of $z$ in the interior, and an eigenfunction 
in the sense that
$$
(\Delta - \zeta(n-\zeta))E_\zeta(\cdot, y') = 0.
$$ 
This follows immediately from $(\Delta - \zeta(n-\zeta))G_\zeta(\cdot, z') =
\delta_{z'}$.  Or one can just check this explicitly, since
$$
E_\zeta(z,y') = c_\zeta \Bigl[\frac{x}{x^2 + |y-y'|^2}\Bigr]^{\zeta}
$$

Given $f\in \cinf_c(\bbR^n)$, one can form a solution to $\Delta u = 
\zeta(n-\zeta)u$ by integrating this generalized eigenfunction
$$
u(z) = 2^{2\zeta}(2\zeta-n) \int_{\bbR^n} E_\zeta(z,y') f(y') \>d^ny.
$$

\begin{proposition}\label{modelp}
For $\re \zeta = n/2$, $\zeta \ne n/2$, $u$ has asymptotic behavior
$$
u(z) = x^{n-\zeta} f(y) + x^{\zeta} f'(y) + O(x^{n/2+1}),
$$
near $x=0$.  Moreover $f' = S_\zeta f$, where $S_\zeta$ is the zeroth order
pseudodifferential operator on $\bbR^n$ with total symbol 
$$
a(y,\xi) = 2^{n-2\zeta} \frac{\Gamma(\frac{n}2-\zeta)}{\Gamma(\zeta-\frac{n}2)}
\>|\xi|^{2\zeta-n}.
$$
\end{proposition}

\begin{proof}
Define $W_x\in \calS'(\bbR^n)$ by 
$$
W_x(|w|) = \frac{x^\zeta}{(|w|^2 + x^2)^\zeta}. 
$$
Then the partial Fourier transform of $u$ is $\hat u(x,\xi) = c_\zeta 
\hat W_x(|\xi|) \hat f(\xi)$.  The distributional Fourier transform 
$\hat W_1(|\xi|)$
is the analytic continuation of a Bessel function, so $\hat W_x$ may be 
analyzed with standard tricks 
(see for example Mandouvalos \cite{Man83}).  
For $\zeta$ as above, as $x\to 0$ we have
$$
\hat W_x(\xi) = \pi^{n/2} \frac{\Gamma(\zeta-\frac{n}2)}{\Gamma(\zeta)} 
\> x^{n-\zeta} + 
\pi^{n/2} \frac{\Gamma(\frac{n}2-\zeta)}{\Gamma(\zeta)} 
\> x^{\zeta} \> \Bigl(\frac{|\xi|}{2}\Bigr)^{2\zeta - n} + O(x|\xi|)
$$
\end{proof}

\section{Stretched products and distributions}\label{spdsec}

The Laplacian $\Delta_g$ for an asymptotically hyperbolic metric may be 
regarded as an operator on $X$ which is elliptic in the interior but 
degenerates uniformly at the boundary.  It belongs to $\diff^m_0$, the 
enveloping algebra of the space of 
smooth vector fields on $X$ which vanish at the boundary.  In local coordinates
$z=(x,y)$, with $x$ a boundary defining function, such operators have the form
$$
\sum_{|\alpha|\le m} a_{\alpha}(z) \bigl(x\deriv{z}\bigr)^\alpha.
$$
In order to even state our main result, we need to review the
calculus of distributions used to analyze the inverses of such operators.

As usual, in order to gave a nice symbol map, we want to consider operators acting
on half-densities.  Over a manifold $W$ let $\hfd(W)$ denote the bundle of
half-densities.  We will simply write $\hfd$ when the manifold is clear from
context. 

\subsection{Stretched product}
Let $\cinfd(X; \hfd)$ be the space of smooth half-densities vanishing to infinite 
order on $\bX$, and $\calD'(X; \hfd)$ the space of distributional half-densities
extendible across $\bX$.  
Linear continuous operators $\cinfd(X; \hfd) \to \calD'(X; \hfd)$ have
Schwartz kernels which are extendible distributional half-densities on $\xx$.   
Note that $\xx$ is a manifold with two boundary hypersurfaces, plus a corner 
$\bxbx$ where they meet.  A crucial element of the parametrix construction of
\cite{MM} is the resolution of singularities of Schwartz kernels at the submanifold
$\bdiag$ where the diagonal $\diag\subset\xx$ meets the corner.  
This is done using the technique of
\cite{Me81}, the introduction of a ``stretched product'' $\sxx$.
$\sxx$ is obtained from $\xx$ by blowing up $\bdiag$, which essentially means
that for each point on $\bdiag$ we keep track of a direction of approach
from $\xx$.  The new manifold has three boundary hypersurfaces: The left and right
faces $F_f$ and $F_r$, corresponding to $\bX\times X$ and $X\times \bX$ in the
original product, and the front face $F_f$, which is the replacement of
$\bdiag$.

Locally this is just the introduction of polar coordinates.  Let $(x,y)$ be
coordinates in $X$, with $x$ a defining function for $\bX$.  Then 
$(x,y,x',y')$ give a set of coordinates for $\xx$ near the boundary of
$\diag$ ($\bdiag$ being given locally by $\{x=x'=0,\; y=y'\}$).  
We introduce polar coordinates around $\bdiag$:
$$
r = \sqrt{x^2 + (x')^2 + |y-y'|^2}, \qquad (\eta,\eta',\theta) =
\frac{(x,x',y-y')}{r}.
$$
Here $(\eta,\eta',\theta)$ lives in a closed quarter sphere of $S^{n+1}$, since
$\eta, \eta'\ge 0$.  The full set of local coordinates for $\sxx$ is
$(r,\eta,\eta',\theta,y)$, and the faces are locally given by
$$
F_l = \{\eta = 0\},\qquad F_r = \{\eta'=0\}, \qquad F_f = \{r=0\}.
$$ 

The global description of the stretched product is as follows.  We take $\xx$,
remove $\bdiag$, and then in its place we glue the doubly inward pointing 
spherical normal bundle of $\bdiag$. This is given a smooth structure using the 
polar coordinate patches introduced above.  Clearly there is a smooth map
$b:\sxx\to\xx$  which collapses the blow up.  In terms of this map, 
$$
F_l = b^{-1}(\bX\times X),\quad F_r = b^{-1}(X\times \bX),\quad F_f = b^{-1}
(\bxbx). 
$$
Let $\phi_l, \phi_r, \phi_f$ be fixed defining functions for the faces 
$F_l, F_r, F_f \subset \sxx$, respectively.

For many purposes it is convenient to use the projective coordinates
$(x,y,t,u)$, where $t = x/x'$ and $u = (y-y')/x'$.  In these coordinates the left,
right, and front faces are $t=0$, $t=\infty$, and $x=0$, respectively.
 
Given a linear continuous operator $A: \cinfd(X; \hfd) \to \calD'(X; \hfd)$, let
$\krnl{A}$ be its Schwartz kernel.  We will characterize such operators by the
pull-back of the kernel to $\sxx$:
$$
\tkrnl{A} := b^*\krnl{A}.
$$
Let $\pi_l$ and $\pi_r$ be the maps $\sxx \to X$ which correspond to 
projections onto the left and right factors in the interior.  
Then given $\tkrnl{A}$, we recover the action of $A$ by
$$
A\mu = \pi_{l*} (\tkrnl{A} \pi_r^* \mu),
$$
for $\mu\in\cinf(X,\hfd)$.
To illustrate this in the local coordinates $(x,y,t,u)$, suppose
$$
\krnl{A} = k(x,y,x',y') \> \Bigl| \frac{dx\>dy\>dx'\>dy'}{(xx')^{n+1}} 
\Bigr|^{1/2}.
$$
Here we include the factors of $x$ and $x'$ in the denominator because that
is the form of the Riemannian half-density $\mu_g$.
Then
$$
\tkrnl{A} = k(x,y,\tfrac{x}{t}, y-\tfrac{xu}{t}) \> \Bigl|
\frac{dx\>dy\>dt\>du}{tx^{n+1}} \Bigr|^{1/2}.
$$
If $\mu = f(x,y) |\frac{dx\>dy}{x^{n+1}}|^{1/2}$ then
$$
\pi_r^* \mu = f(\tfrac{x}{t}, y-\tfrac{xu}{t}) 
\Bigl|\frac{dt\>du}{t}\Bigr|^{1/2}.
$$
The product $\tkrnl{A} \pi_r^* \mu$ is a density in the $t$ and $u$ variables
and can be pushed forward to give
\begin{equation}\label{tkact}
A\mu = \int k(x,y,\tfrac{x}{t}, y-\tfrac{xu}{t}) f(\tfrac{x}{t}, 
y-\tfrac{xu}{t})
\frac{dt\>du}{t} \;\cdot \Bigl|\frac{dx\>dy}{x^{n+1}}\Bigr|^{1/2}.
\end{equation}

Finally we note that for the identity operator $I$,
\begin{equation}\label{ident}
\tkrnl{I} = \delta(t-1) \delta(u) \>
\Bigl|\frac{dx\>dy\>dt\>du}{x^{n+1}}\Bigr|^{1/2}
\end{equation}

\subsection{Distributions}
The Riemannian half-density $\mu_g$ is a section of the singular density
bundle
$(\rho\rho')^{-(n+1)/2}\hfd(\xx)$.  The pull-back of this bundle to $\sxx$ is
$$
b^*\Bigl[ (\rho\rho')^{-(n+1)/2}\hfd(\xx)\Bigr] = \tau \hfd(\sxx),
$$
where
$$
\tau = (\phi_l \phi_r \phi_f)^{-(n+1)/2}.
$$
Accordingly, we will define spaces of operators whose kernels are sections 
of this bundle.

Let $\dlift$ denote the closure of the lift to $\sxx$ of the interior of the
diagonal in $\xx$.  Note that $\dlift$ intersects the front face 
transversally.  
We define $\cnrml{m}$ to be the space of distributional sections of
$\tau\hfd(\sxx)$ which are conormal to $\dlift$ of degree $m$ as in \cite{H}
and which vanish to all orders at $F_l$ and $F_r$. Define
$$
\op{m}{} := \{A: \cinfd(X, \hfd) \to \calD'(X, \hfd);\; \tkrnl{A} \in 
\cnrml{m} \}.
$$
Because of the vanishing at the left and right faces, the factors of 
$\phi_l$ and $\phi_r$ in $\tau\hfd$ are irrelevant to the definition.
To see the significance of the $\phi_f$, note from
from (\ref{ident}) that $\tkrnl{I}$ is a smooth section of
$\tau\hfd$ but would be singular as a section of $\hfd$ (in those
coordinates $\phi_f = x$).  So the definition above gives the desired
fact that $I\in \op{0}{}$.

The topology of $\cnrml{m}$ is defined as follows.  
Since the distributions are extendible it is convenient to 
consider $(\sxx)^2$, which is the double of $\sxx$ across this face.
Near the diagonal of $(\sxx)^2$ we take local Fourier transforms
and use the topology of the standard symbol spaces $S^m$.  
Away from the diagonal and the topology is that of 
$(\tilde\phi_l \tilde\phi_r)^\infty \cinf((\sxx)^2)$, where 
$\tilde\phi_l$ and $\tilde\phi_r$ are any smooth extensions of 
$\phi_l$ and $\phi_r$.  The topology of $\cnrml{m}$ is simply
the restriction of this topology on $(\sxx)^2$.
This definition makes $\cnrml{m}$ a complete locally
convex vector space.

It is important to note that operators in $\op{-\infty}{}$ have kernels 
which are smooth when
pulled back to $\sxx$, but not necessarily smooth on $\xx$.

Our constructions will also involve kernels with polyhomogeneous conormal
singularities at the boundaries.  We'll review the facts we need;
see \cite{Me92}, \cite{MeAPS}, \cite{MeBk} for more complete expositions.  
Let $W$ be a manifold with corners.  Label the boundary faces
$1,\dots,k$, with corresponding boundary defining functions $\phi_j$.  
Polyhomogeneous conormal singularities are described by the powers of 
$\phi_j$ and $\log \phi_j$ which occur in the expansions at each face.  For 
each $j = 1,\ldots, k$ an {\it index set} $E_j$, a countable discrete subset 
of $\bbC\times \bbN_0$.  The collection $\calE = \{E_1,\dots, E_k\}$ is called
an {\it index family} for $W$.	 The space of {\it polyhomogeneous conormal
distributions}\footnote{We omit the usual designation ``phg'' from the notation
because all distributions in this paper are polyhomogeneous at boundaries.} 
$\phgd{\calE}(W)$ consists of functions $u$ which are smooth on the interior 
of $W$ and which near the $j-$th boundary face have an asymptotic expansion 
of the form:
\begin{equation}\label{phgexp}
u \sim \sum_{(a,l)\in E_j} \sum_{0=1,\dots,l} \phi_j^a (\log \phi_j)^l u_{a,l},
\end{equation}
where the $u_{a,l}$ are smooth.  (To give a proper definition one must take 
somewhat more care with the corners; we refer the reader to \cite{MeBk}.)
To insure that this expansion makes sense,
for each index set $E$ it is required that the set
$$
E \cap (\{\re a < M\}\times \bbN_0)
$$
is finite for all $M\in\bbR$.  If $u$ vanishes to infinite order at the $j-$th
face, then we write $E_j = \infty$.  
The definitions are local and can be applied to sections of line
bundles.

For simplicity, we will abbreviate $E = \{(a,0)\}$ simply as $a$.
This indicates singularities of a very simple form:
$$
\phgd{a_1,\dots,a_k}(W) = \phi_1^{a_1} \dots \phi_k^{a_k} \>\cinf(W).
$$

The spaces of polyhomogeneous conormal distributions are given 
topologies as follows.
We first define the $C^\infty$ seminorms.  Let $\calV_b(W)$ be the set of smooth
vector fields tangent to all boundary faces, and choose a set $\{V_1,\dots,V_l: \;
V_j\in \calV_b(W)\}$ which together span $\calV_b(W)$ everywhere. 
The index sets have a partial ordering:
$(a,l) < (b,k)$ if $\re a < \re b$ or if $a=b$ and $l>k$.  
So given an index family $\calE$ we define $\Emin$ by choosing a smallest 
member of each index set (this may not be unique).  
For $k\in \bbN_0$ we define the norm on $\phgd{\calE}(W)$:
$$
\norm{u}_{k;\Emin} =  \sum_{|\alpha| \le k} \sup_\sxx \bigl|
V^\alpha\>\phi^{-\Emin} u \bigr|,
$$
where 
$$
\phi^{-\Emin} := \prod_j \phi_j^{-e_j} (\log\phi_j)^{-m_j},
$$
for $\Emin = \{(e_1,m_1), (e_2,m_2), \dots\}$.
The imaginary parts of the $e_j$ are irrelevant in this definition, 
so the non-uniqueness in the choice of $\Emin$ doesn't matter.

The full topology is defined inductively, using $\cinf$ seminorms on the
coefficients $u_{a,l}$ from the expansion (\ref{phgexp}) in local coordinates,
and then seminorms of the form $\norm{\cdot}_{k;\Emin}$ on the remainders when
leading terms in the expansion have been subtracted off.  In essence, this
topology controls all coefficients and all remainders from the asymptotic
expansions.  With this topology $\phgd{\calE}(W)$ is also a complete locally
convex vector space \cite{MeBk}.

Returning now to the lifts of kernels from $\xx$ to $\sxx$, let
$\calE$ be an index family for $\sxx$, with the boundary faces ordered
left, right, front.  We define
$$
\op{}{\calE} := \{A: \cinfd(X,\hfd) \to \calD'(X,\hfd);\; \tkrnl{A} \in 
\phgk{\calE}\}.
$$
Most commonly, $\calE$ will be of the form $\{a,b,c\}$.  We also denote
$$
\op{m}{\calE} = \op{m}{} + \op{}{\calE},
$$
with the topology inherited from the two pieces.

Another space we will need consists of operators with smooth kernels on $\xx$,
which are polyhomogeneous conormal at the left and right boundaries.  (Such 
operators have nothing to do with the stretched product.)  Given an index 
family $\calI$ for $\xx$, define
\begin{equation*}
\begin{split}
\vrop{\calI} &:= \{A: \cinfd(X,\hfd) \to \calD'(X,\hfd);\;\\
&\hskip1in \krnl{A} \in \phgd{\calI}(\xx; (\rho\rho')^{-(n+1)/2}\hfd) \}.
\end{split}
\end{equation*}
The normalization of the half-density bundles is such that
$$
\vrop{a,b} \subset \op{}{a,b,a+b}.
$$

The composition of operators in $\op{}{\calE}$ combines index families in a
straightforward but non-trivial way.  Given index sets $E_1, E_2$, the sum $E_1
+ E_2$ has the obvious meaning:
$$
E_1 + E_2 = \{(a+b,k+l):\; (a,k) \in E_1,\> (b,l) \in E_2) \}.
$$
But because compositions can introduce extra logarithmic singularities, 
an extended notion of union is required.  If $b\notin a+\bbZ$, then we simply set
$$
\{(a,k)\} \xU \{(b,l)\} := \{(a,k), (b,l)\}.
$$
But if $b \in a+\bbN_0$, then we define
$$
\{(a,k)\} \xU \{(b,l)\} :=  \{(a,k), (b,k+l+1)\}.
$$ 
This notion is extended to full index
sets in the obvious way. We also define
$$
\re E := \min\{\re a: (a,k) \in E\}.
$$
The index set $E = \infty$ has the additive property suggested by the 
notation, but  behaves as the empty set in unions:
$$
\infty + F = \infty, \qquad \infty\xU F = F.
$$

The following result is taken from Theorems 3.15
and 3.18 of Mazzeo's paper (we use a slightly different convention 
for the index sets).
\begin{theorem}\label{compthm}\cite{Ma91}\footnote{The continuity of the 
compositions was not actually considered explicitly in \cite{Ma91}.  The techniques
used in \cite{Ma91} are shown to be continuous in \cite{MeBk}.}
For $m, k\in \bbZ$ and index sets given by $\calE = \{E_1, E_2, E_3\}$, $\calF =
\{F_1,F_2,F_3\}$ such that $\re E_2 + \re F_1 > n$, composition gives a 
continuous map
$$
\op{m}{\calE} \times \op{k}{\calF} \to \op{m+k}{\calE\circ\calF},
$$
where
$$
\calE\circ\calF := \{E_1 \xU (F_1 + E_3),\; F_2 \xU (E_2 + F_3),\; (E_3 + F_3)
\xU (E_1 + F_2)\}.
$$

Also, for $m\in \bbZ$ and index sets $\calE = \{E_1, E_2, E_3\}, \calI =
\{I_1,I_2\}$ such that $\re E_2 + \re I_1 > n$, composition gives a continuous map
$$
\op{m}{\calE} \times \vrop{\calI} \to \vrop{\calE\circ\calI},
$$
where 
$$
\calE\circ\calI := \{E_1 \xU (I_1 + E_3), I_2\}.
$$
\end{theorem}

To compose operators, the respective kernels are pulled back to $X\times X\times X$,
multiplied together, and then pushed forward.  To keep track of the conormal singularities
in the first composition formula above, 
$X\times X\times X$ is replaced by a blown up version, denoted by $\sxxx$.  
This space is equipped with maps $\beta_{ij}: \sxxx\to \sxx$, $i,j = 1,2,3$, which
correspond on the interior to projections onto the $i$-th and $j$-th factors.  The crucial
feature is that these maps are boundary fibrations in the sense of
\cite{MeBk}.    The composition of $f,g \in \phgk{*}$ is given by
\begin{equation}\label{lkcomp}
f\circ g := \beta_{13*}  [\beta_{12}^* f \cdot \beta_{23}^* g].
\end{equation}

One can also study the action of $\op{m}{\calE}$ on conormal functions on $X$
by similar methods.  We will need only the following.

\begin{proposition}\label{opact}\cite{Ma91}
For $A\in \op{m}{\calE}$ and $u\in \dot\cinf(X; \hfd)$, $Au$ is well-defined
and
$$
Au \in \phgd{E_1}(X;  \rho^{-(n+1)/2}\hfd),
$$
where $E_1$ is the index set in $\calE$ corresponding to the left face.
\end{proposition}

Given the metric $g$ the natural Hilbert space is 
$L^2(X, \mu_g^2)$.  The map $f\mapsto f \mu_g$ removes the dependency 
on the metric and 
defines an isometry $L^2(X,\text{vol}_g) \to \hfdhil$.  
We will need to consider weighted $L^2$ spaces of the form
$\rho^\delta \hfdhil$ for $\delta>0$.

\begin{proposition}\cite{Ma91}
Operators in $\op{}{\calE}$, where $\calE = \{E_1, E_2, E_3\}$, 
are compact on $\rho^\delta \hfdhil$ provided that 
$\re E_1 > \frac{n}{2} + \delta$, 
$\re E_2 > \frac{n}{2} - \delta$, and $\re E_3 \ge 0$.
\end{proposition}

\subsection{Symbol map}
We come now the main reason for the introduction of the stretched product.
Consider the operators $x\del_{z_j}$ in local coordinates $z = (x,y)$. 
Let $A\in \diff^1_0(X; \hfd)$ be given locally by
$$
A: f(z) \>\Bigl|\frac{dx\>dy}{x^{n+1}}\Bigr|^{1/2} \mapsto x \deriv{z_j}f(z)\> 
\Bigl|\frac{dx\>dy}{x^{n+1}}\Bigr|^{1/2}.
$$
From (\ref{tkact}) we see that
$$
\tkrnl{A} = t \deriv{v_j} [\delta(t-1) \delta(u)] \>
\Bigl|\frac{dx\>dy\>dt\>du}{tx^{n+1}}\Bigr|^{1/2},
$$
where $v = (t,u)$.  Note that $\tkrnl{A}$ no longer degenerates; it is
transversally elliptic to $\dlift$, all the way down to the front face
$\{x=0\}$.

Given a metric $g$ and $\zeta\in\bbC$, we define the operator
$P_\zeta \in \diff^2_0(X; \hfd)$ by
$$
P_\zeta (f \mu_g) := \bigl[ (\Delta_g - \zeta(n-\zeta))\>f \bigr] \> \mu_g,
$$
where $f\in \cinf(X)$ and $\mu_g$ is the Riemannian half-density associated to
$g$ (note that $\mu_g$ is singular at the boundary).

\begin{proposition}\label{ellip}\cite{MM}
$\tkrnl{P_\zeta}$ is uniformly transversally elliptic to $\dlift$.   
\end{proposition}

\begin{proof}
This may be checked in local coordinates $z=(x,y)$ as above, with $x=\rho$. 
Locally $\Delta_g$ takes the form
\begin{equation}\label{Laploc}
\Delta_g = - \sum_{j,k} \frac{1}{\sqrt{h}} \bigl(x\deriv{z_j}\bigr)
h^{jk} \sqrt{h} \bigl(x\deriv{z_k}\bigr) + n \sum_k h^{0k}
\bigl(x\deriv{z_k}\bigr),
\end{equation}
where $h_{ij} = x^2 g_{ij}(x,y)$.
If $\mu_g = w(x,y) |\frac{dx\>dy}{x^{n+1}}|^{1/2}$, then 
$$
\tkrnl{\Delta_g} = L[\delta(t-1)\delta(u)] \frac{w(x,y)}{w(\frac{x}{t}, 
y - \frac{tu}{x})} \> \Bigl|\frac{dx\>dy\>dt\>du}{tx^{n+1}}\Bigr|^{1/2}.
$$
where $L$ is the differential operator:
\begin{multline}
L = - \sum_{j,k} h^{jk} \bigl(t\deriv{v_j}\bigr)
\bigl(t\deriv{v_k}\bigr) - \sum_{j,k} \frac{x}{\sqrt{h}}
\deriv{z_j}(h^{jk} \sqrt{h}) \bigl(t\deriv{v_k}\bigr) \\
+ n \sum_k h^{0k} \bigl(t\deriv{v_k}\bigr) ,
\end{multline}
with $v = (t,u)$ and $h = h(x,y)$.  
Now as the front face $\{x=0\}$ is approached, the principal
term is $t^2 \sum h^{jk} \deriv{v_j}\deriv{v_k} \delta(t-1) \delta(u)$. 
This proves the proposition, since $t=1$ on $\dlift$ and $h$ is smoothly
extendible across $x=0$.
\end{proof}

Proposition \ref{ellip} shows that in local coordinates $\tkrnl{\Delta_g}$ 
looks like the kernel of an elliptic operator which does not degenerate at the
boundary.  In fact there is a symbol map  for elements of $\op{m}{}$, given by
taking the symbol of the lifted kernel as a conormal distribution on $\sxx$.
Proposition \ref{ellip} is then simply the statement that the 
$\sigma(\Delta_g)$ defined in this way is invertible.

The symbol of a conormal distribution is a half-density on the conormal bundle
of the singular set (see \cite{H}).  The conormal bundle of $\dlift$ in $\sxx$
is naturally identified with a bundle
$\tilde T^*X$ over $X$, the so-called {\it compressed cotangent bundle}. 
Locally this bundle is spanned by differentials of the form $dx/x, dy_j/x$. 
$\tilde T^*X$ carries a natural symplectic form, and thus a natural density.  
So symbols can be invariantly identified with functions on $\tilde T^*X$ (which
is why we use operators on half-densities).  These functions will lie in the
standard symbol spaces $S^m(\tilde T^*X)$.  

\begin{proposition}\cite{MM}\label{symseq}
The symbol map gives an exact sequence:
$$
0\longrightarrow \op{m-1}{} \longrightarrow \op{m}{}
\xrightarrow{\>\sigma\>} S^m(\tilde T^*X)/S^{m-1}(\tilde T^*X)
\longrightarrow 0,
$$
and the symbol is multiplicative:
$$
\sigma(AB) = \sigma(A) \cdot \sigma(B).
$$
\end{proposition}

Using Propositions \ref{ellip} and \ref{symseq}, we may invert
$P_\zeta$ symbolically, with an error in $\op{-\infty}{}$.  
As noted above, the kernel of such an operator is in general 
singular at the corner of $\xx$.
So we need a way to refine the parametrix near the front face.

\subsection{Normal operator}
The tool used to accomplish this problem is the normal
operator at the front face. Given a point $p\in\bX$, let $T^+_pX\cong
\bbR^{n+1}_+$ be the (closed) inward half of $T_pX$.   
Given $A\in \op{m}{\calE}$, $p\in \bX$, the restriction of $\tkrnl A$ to 
$b^{-1}(p,p) \subset F_f$ may be regarded as an operator on $\hfd(T^+_pX)$
by convolution.  This operator is denoted by $N_p(A)$.

It is easiest to define $N_p(A)$ in local coordinates.  
Let $(x,y,t,u)$ be local coordinates for $\sxx$, with $p$ given by $(0,y_0)$.  
We also use $(x,y)$ as coordinates for $T^+_pX$.
Then if 
$$
\tkrnl{A} = k(x,y,t,u) \> \Bigl|\frac{dx\>dy\>dt\>du}{tx^{n+1}}\Bigr|^{1/2},
$$
the action of the normal operator is
$$
N_p(A): f(x,y) \> \Bigl|\frac{dx\>dy}{x^{n+1}}\Bigr|^{1/2}
\mapsto \int k(0,y_0,t,u) f(\frac{x}{t}, y-\frac{xu}{t}) \frac{dt\>du}{t}
\;\cdot \> \Bigl|\frac{dx\>dy}{x^{n+1}}\Bigr|^{1/2},
$$
for $f\in \cinf(T^+_pX)$.  Note that 
$$
N_p(I) = I.
$$

For differential operators in $\diff^*_0(X)$, the definition 
amounts to ``freezing coefficients'' at the boundary:
if $A$ is given by $Af \mu \mapsto a(x,y) (x\del_z)^\alpha f \>\mu$
for some non-vanishing $\mu\in \hfd(X)$, then we have
\begin{equation}\label{Ndiff}
N_p(A): f \> \Bigl|\frac{dx\>dy}{x^{n+1}}\Bigr|^{1/2}
\mapsto  a(0,y_0) (x\del_z)^\alpha f \> \Bigl|\frac{dx\>dy}{x^{n+1}}\Bigr|^{1/2}
\end{equation}
Let $g_p$ be the (hyperbolic) metric on $T^+_pX$ given by $x^{-2} h|_p$, 
where $h|_p$ is interpreted as a constant matrix.  It follows from
(\ref{Laploc}) and (\ref{Ndiff}) that
$$
N_p(P_\zeta) = \Delta_{g_p} - \zeta(n-\zeta).
$$

As a function of $p$, $N_p(A)$ is clearly smooth in the interior of $F_f$.  In
fact, if we think of $N(A)$ as a function on $F_f$, then $N$ maps $\op{}{a,b,0}$
to $\phgd{a,b}(F_f)$.  Clearly $N$ is continuous in the topology defined above
for polyhomogeneous distributions, since it just amounts to reading off the
leading asymptotic coefficient at the front face.

\begin{proposition}\label{nrmlprop}\cite{MM}
The normal operator gives an exact sequence
$$
0\longrightarrow \op{}{a,b,1} \longrightarrow \op{}{a,b,0} \xrightarrow{\>N\>}
\phgd{a,b}(F_f) \longrightarrow 0.
$$
For $P\in \diff_0^*(X, \hfd)$, 
$$
N_p(P\cdot K) = N_p(P)\circ N_p(K).
$$
\end{proposition}

\section{Resolvent and scattering operator}\label{rsosec}
As in \S\ref{spdsec}, let $P_\zeta$ be the operator $\Delta_g - \zeta(n-\zeta)$,
acting on half-densities through the Riemannian half-density.
The resolvent, $R_\zeta := P_\zeta^{-1}$ exists and is bounded for $\re\zeta$
sufficiently large, because $\Delta_g$ is self-adjoint and positive.

\begin{theorem}\label{mmthm} \cite{MM}, \cite{Ma91}
$R_\zeta$ may be analytically continued to a meromorphic family on the
domain $\zeta \in \bbC\backslash\{\frac12(n-\bbN)\}$ with
$$
R_\zeta \in \op{-2}{\zeta,\zeta,0} + \vrop{\zeta,\zeta}.
$$
\end{theorem}

The method of proof is as follows.  For $\zeta \notin \frac12(n-\bbN)$ a
parametrix 
$$
M_\zeta \in \op{-2}{\zeta,\zeta,0}
$$ 
is constructed, such that $P_\zeta M_\zeta - I = E_\zeta \in
\vrop{\infty,\zeta}$.  This error term has two crucial features. 
First of all, it's compact on $\rho^\delta L^2(X, \rho^{-(n+1)/2} \hfd)$ 
for $\delta>0$, so $I + E_\zeta$ can be inverted meromorphically by analytic 
Fredholm theory.  Secondly, if we
define $F_\zeta$ by setting $I + F_\zeta = (I+E_\zeta)^{-1}$, then $F_\zeta$ is
also in $\vrop{\infty,\zeta}$.  This in turn implies that $R_\zeta =
M_\zeta  + M_\zeta F_\zeta \in \op{-2}{\zeta,\zeta,0} + \vrop{\zeta,\zeta}$.

We turn next to some facts concerning solutions of
$\Delta_g u = \zeta(n-\zeta) u$.  The first is essentially a uniqueness 
result. 
Although well-known, it doesn't seem to be proven in the literature so we will
give a proof.  Let $\Spec(\Delta_g)\subset (0,n^2/4)$ be the point spectrum of 
the Laplacian.

\begin{proposition}\label{uniqprop}
Assume $\re\zeta \ge n/2$, $\zeta\ne n/2$ and $\zeta(n-\zeta)
\notin\Spec(\Delta_g)$.   Suppose that the function
$u$ that solves $\Delta_g u = \zeta(n-\zeta)u$ and has asymptotic behavior
\begin{equation}\label{uniq}
u = \rho^\zeta f + O(\rho^{\zeta+\epsilon}),
\end{equation}
for $f\in\cinf(\bX)$, with $\epsilon>0$.  Then $u=0$.
\end{proposition}

\begin{proof}
For $\re\zeta>n/2$ the expansion (\ref{uniq}) implies that $u$ is in
$L^2(X,\text{vol}_g)$ so in this case the statement is tautological.

For the rest of the proof we assume $\re\zeta=n/2$.
We will show that $\Delta_g u = \zeta(n-\zeta)u$ 
and (\ref{uniq}) together imply $f=0$.  
This will imply $u\in L^2(X, \text{vol}_g)$, and we conclude $u=0$ as before.

We use a version of the boundary pairing argument of \cite{Me94}.  Suppose that 
$\Delta_g u = \zeta(n-\zeta)u$.  Let $\psi\in \cinf(\bbR_+)$ be a  
cutoff function so that 
$\psi(t) = 0$ for $t\le 1$ and $\psi(t) = 1$ for $t\ge 2$.  
By the self-adjointness of $\Delta_g - \zeta(n-\zeta)$, we have
$$
\int_X \Bigl([\Delta_g,\psi(\lambda\rho)] u\Bigr)\> \overline{u} \>dg = 0.
$$
So in particular
\begin{equation}\label{eplim}
\lim_{\lambda\to\infty} \int_X \Bigl([\Delta_g,\psi(\lambda\rho)] u\Bigr)
\> \overline{u} \>dg = 0.
\end{equation}
In local coordinates $(x,y)$ with $x=\rho$,
$$
\Delta_g = - (x\del_x)^2 + nx\del_x + p(x,y,x\del_y) + x T,
$$
where $T\in \diff_0^2(X)$.
Thus if $u$ satisfies (\ref{uniq})
$$
[\Delta_g,\psi(\lambda x)]u = \Bigl[(n-2\zeta-1) \lambda
\psi'(\lambda x)  - \lambda^2 x \psi''(\lambda x) \Bigr] x^{\zeta+1} f
+ \lambda \psi'(\lambda x) \times O(x^{n/2+1+\epsilon}).
$$
Substitute this expression back into (\ref{eplim}).   After we perform the $x$
integration, only the term involving $f$ will survive as 
$\lambda\to\infty$ (recall that $dg = dh/x^{n+1}$).  The conclusion is that
$$
0 = (n-2\zeta) \int_{\bX} |f|^2 \>dh|_\bX,
$$
so $f$ must vanish.
\end{proof}
The second result describes the asymptotic behavior of solutions of 
$\Delta_g u = \zeta(n-\zeta) u$ and allows us to define 
the scattering operator.  Once again this is well-known, and more
general results on asymptotic expansions of generalized eigenfunctions
may be found in \cite{Ma91}.

\begin{proposition}\label{aeprop}
Let $\re\zeta = n/2, \im\zeta\ne 0$. Given $f\in\cinf(\bX)$ there is a unique
solution of $\Delta_g u = \zeta(n-\zeta) u$ which near $\bX$ has the asymptotic
form:
$$
u(z) = \rho^{n-\zeta}f + \rho^{\zeta}f' + O(\rho^{n/2+\epsilon}).
$$
\end{proposition}

\begin{proof}
As in the model case, we define 
\begin{equation}\label{ezdef}
E_\zeta(q,p') = \lim_{q'\to p'} \rho(q')^{-\zeta} 
\frac{\krnl{R_\zeta}(q,q')}{\mu_g \mu'_g},
\end{equation}
and note that
$$
u(q) = 2^{2\zeta}(2\zeta-n) \int E_\zeta(q,p') f(p') \>dh|_\bX
$$
is a solution of $\Delta_g u = \zeta(n-\zeta) u$.  
Choose local coordinates $(x,y)$ around $p = (0,y_0)$ for which $g_p$ is the
standard hyperbolic metric.   Then because 
$N_p(\Delta_g) = \Delta_{g_p}$ and by the exact sequence of Proposition
\ref{nrmlprop}, in the coordinates $(\eta,\eta',\theta,r,y_0)$ we have
\begin{equation}\label{rlocal}
\tkrnl{R_\zeta}(\eta,\eta',\theta,r,y_0) = \Bigl(G_\zeta(\eta,\eta',\theta) +
r(\eta\eta')^\zeta F(\eta,\eta',\theta,r,y_0)\Bigr) \; 
\Bigl|\frac{d\eta\>d\eta'\>d\theta\>dr\>dy}{(\eta\eta'r)^{n+1}}\Bigr|^{1/2},
\end{equation}
where $F$ is smooth.  Here $G_\zeta$ is the model resolvent, which depends only
on $\eta,\eta'$, and $\theta$ because in these coordinates the hyperbolic 
distance $d$ is given by
$$
\cosh d = 1 + \frac{|\eta-\eta'|^2 + \theta^2}{2\eta\eta'}.
$$
From (\ref{rlocal}) and (\ref{ezdef}) we compute that
\begin{equation}\label{ezlead}
E_\zeta(x,y_0,y') = \frac{c_\zeta x^\zeta}{r^{2\zeta}}
 + \frac{x^\zeta}{r^{2\zeta-1}} F(x/r,0,(y_0-y')/r,r,y_0),
\end{equation}
where $r = \sqrt{x^2 + (y_0-y')^2}$ and $c_\zeta$ is the constant
(\ref{czeta}).
The analysis proceeds as in Proposition \ref{modelp}.  

Uniqueness of the solution follows immediately from Proposition \ref{uniqprop}.
\end{proof}

Using Proposition \ref{aeprop}, for $\re\zeta = n/2$
we define $S_\zeta$ to be the operator which
maps $f\mapsto f'$.  From the proof of the proposition 
it is clear that the Schwartz kernel with respect to the Riemannian
density on $\bX$ induced by $h$ is 
\begin{equation}\label{sckdef}
K_{S_\zeta}(p,p') = 2^{2\zeta} (2\zeta-n) \lim_{z\to p} \lim_{z'\to p'} 
\rho(z)^{-\zeta} \frac{\krnl{R_\zeta}(z,z')}{\mu_g \mu'_g} \rho(z')^{-\zeta},
\end{equation}
for $p\ne p' \in \bX$.  We can adopt (\ref{sckdef}) more generally 
to define $S_\zeta$ meromorphically in $\zeta$.\footnote{Note that this
limiting process will introduce new poles in $S_\zeta$ which were not poles of
$R_\zeta$, because $\krnl{S_\zeta}$ must be interpreted as a distribution.  This issue
of resolvent vs. scattering poles will be dealt with in \cite{BP}.}  
Note that Proposition \ref{aeprop} implies that for
$\re\zeta = n/2$, 
$$
S_{n-\zeta}S_\zeta = I.
$$
Since $S(\zeta)$ is a meromorphic family, this relation continuous to hold
for all $\zeta\notin \frac12(n-\bbN)$.  

From (\ref{sckdef}) and the form of $\tkrnl{R_\zeta}$ we see that
in local coordinates $(y,y')$ the kernel of the scattering operator has the form
\begin{equation}\label{localks}
\krnl{S_\zeta} = r^{-2\zeta} F(r,\theta,y) + G(y,y'),
\end{equation}
where $r=|y-y'|$, $\theta = (y-y')/r$, and $F$ and $G$ are $\cinf$ in their respective
variables.   Let $\psbop{a}$ be the set of one-step pseudifferential operators of order
$a \in \bbC$.  That is, local Fourier transforms near the diagonal give symbols with
full asymptotic expansions of the form
$$
a(x,\xi) \sim \sum_{j=0}^\infty a_j(x,\xi) |\xi|^{a-j},
$$
where the $a_j$'s are all homogeneous of degree zero.
The topology on $\psbop{a}$ is given by 
applying the seminorms from the symbol class $S^0$ to each coefficient $a_j$ in the
local expansions near the diagonal (together with $\cinf$ seminorms away from the
diagonal).

\begin{proposition}
$$
S_\zeta \in \psbop{2\zeta-n}.
$$
The principal symbol of $S_\zeta$ is 
$(|\xi|_{h|\bX})^{2\zeta-n}$ times a function meromorphic in $\zeta$ and
independent of $g$.
\end{proposition}

\begin{proof}
The first statement follows immediately from the local form (\ref{localks}).
The existence of the one-step expansion corresponds to the smoothness of 
$F$ as a function of $r$. 
The principal symbol of $S_\zeta$ is derived directly from the limiting 
form for $E_\zeta$ given in local coordinates by (\ref{ezlead}).
\end{proof}

There are standard relations between $R_\zeta$, $E_\zeta$, and $S_\zeta$.  For
proofs of the following, see for example Theorems 5.3 and 6.3 of \cite{P89}.
\begin{proposition}\label{rseprop}
$$
\krnl{R_{n-\zeta}}(q,q') = \krnl{R_\zeta}(q,q') + (n-2\zeta)\mu_g\mu'_g \int_{\bX}
E_\zeta(q, p) E_{n-\zeta}(q',p) \>dh|_{\bX},
$$
and
$$
E_{n-\zeta}(q,\cdot) = S_{n-\zeta}E_\zeta(q,\cdot).
$$
\end{proposition}

\section{Continuity}\label{contsec}

The construction of the parametrix $M_\zeta$ in \cite{MM}
involves summing three asymptotic series, at the lifted diagonal $\dlift$, the front
face, and the left face.   These summations could each be performed so as to insure
continuity in
$g$, but it turns out that we can use an abbreviated construction with only the
first summation.  The reason is that in
the final stage of our argument we will need to use a uniqueness property of the
resolvent to extend continuity from the parametrix to the resolvent.
The uniqueness, which comes from Proposition \ref{uniqprop}, is strong enough to apply
even with a cruder parametrix.  The error term for our parametrix will lie in
$\op{}{\zeta+1,\zeta,1}$ rather than $\vrop{\infty,\zeta}$.  

Let $\metr$ denote the space of smooth, asymptotically hyperbolic metrics on $X$.
The $\cinf$ topology on $\metr$ is defined by seminorms of the form
$$
\norm{\omega} = \sup_X |V_1\dots V_m \rho^2 \omega(Y,Z)|,
$$
where $V_1,\dots,V_m,Y,Z \in \calV(X)$.

As a preliminary, we note the following, which follows easily from the definitions
of the topologies.
\begin{lemma}\label{pzcont}
For any $\zeta\in\bbC$, the map 
$$
\metr \ni g\mapsto P_\zeta \in \op{2}{}
$$
is continuous.
\end{lemma}

Our main goal in this section is to extend this continuity from $P_\zeta$
to its inverse.  

\subsection{The half-plane $\re\zeta\ge n/2$}
\begin{theorem}\label{resthm}
Assume $\re\zeta \ge \frac{n}2$, $\zeta \ne n/2$.  The map 
$$
\metr \ni g\mapsto R_\zeta \in \op{-2}{\zeta,\zeta,0} + \vrop{\zeta,\zeta}
$$
is continuous, except at points where $\zeta(n-\zeta) \in
\Spec(\Delta_g)$.
\end{theorem}

We begin with an elementary topological lemma.
\begin{lemma}\label{approx}
Let $W, d$ be a separable metric vector space with a dense vector subspace
$W_0$.  Suppose $F$ is a continuous map from a topological space
$Y \to W$. Then for any $\epsilon>0$ there exists a continuous map 
$G:Y \to W_0$ such that $d(F(y), G(y)) < \epsilon$ for all $y\in Y$. 
\end{lemma}

\begin{proof}
Since $W$ is separable and $W_0$ is dense, we may choose a countable set 
$w_i \in W_0$ such that $W = \cup_i B_\epsilon(w_j)$, where $B_\epsilon(w_j) = 
\{w\in W: d(w,w_j)<\epsilon\}$.  Let $\phi_j$ be a positive continuous function
such that $0\le \phi_j(w) \le \min\{2^{-j}, 2^{-j}/d(0,w_j)\}$ and $\supp\phi_j
= \overline{B_\epsilon(w_j)}$.  Then $\phi = \sum_j\phi_j$ converges uniformly
on $W$ and hence is continuous.  Also $\phi(w) \ne 0$ for all $w\in W$. Note
that $\sum_j \phi_j(F(y)) w_j$ converges uniformly on $Y$ and so defines a
continuous function $Y\to W$.  Then $G(y) = \frac{1}\phi\sum_j \phi_j(F(y))
w_j$ has the desired properties.
\end{proof}

\begin{proof}[Proof of Theorem \ref{resthm}]
For notational convenience, choose a metric $\cmet_m$ for each $\op{m}{}$,
such that $\cmet_{m} \le \cmet_{m-1}$.

The first stage of the construction is to remove the conormal singularity at $\dlift$;
this is just a standard parametrix construction using the symbol map.
To start, let $Q_{0,0} = I$.  At each inductive step, we are given $Q_{0,j}\in
\op{-j}{}$  which is continuous in $g$ with respect to $d_{-j}$.  
We choose $A_j$ with symbol
$$
\sigma(A_j) = \sigma(Q_{0,j})/\sigma(P_\zeta).
$$
Clearly we may do this so that $A_j$ is a continuous function of $g$ in 
$\op{-2-j}{}$.  Then, using Lemma
\ref{approx}, we may find $E_j \in \op{-\infty}{}$ such that 
$$
\cmet_{-2-j}(\tkrnl{A_j},\tkrnl{E_j}) < 2^{-j}
$$
and $E_j$ is continuous in as an element of $\op{-2-j}{}$.  We take
$$
Q_{0,j+1} = Q_{0,j} - P_\zeta(A_j-E_j) \in \op{-j-1}{},
$$
also continuous, and proceed to the next step.  

Finally, we set
$$
A = \sum_{j=0}^\infty (A_j - E_j),
$$
which converges uniformly in $\op{-2}{}$ and so is continuous as a function of $g$.
Let 
$$
Q_1 = I - P_\zeta A \in \op{-\infty}{}
$$
For any $N$ we can write
$$
Q_1 = Q_{0,N} - \sum_{j=N+1}^\infty P_\zeta (A_j - E_j),
$$
which shows that $Q_1$ is continuous with respect to $d_{-N}$.  Thus $Q_1$ is a continuous
function of $g$ in the topology of $\op{-\infty}{}$.

The next stage is to use the model resolvent on $T^+_pX \cong
\bbH^{n+1}$ to find $B \in \op{}{\zeta,\zeta,0}$ which solves
\begin{equation}\label{istepa}
[\Delta_{g_p} - \zeta(n-\zeta)]\cdot N_p(B) = N_p(Q_1),
\end{equation}
for each $p\in F_f$.  
By the exact sequence of Proposition \ref{nrmlprop}, we then have
\begin{equation}\label{istepb}
P_\zeta B + Q_1 = Q_2 \in \op{}{\zeta,\zeta,1}.
\end{equation}
Note however, that for any $f\in\cinf(X)$,
$$
(\Delta_g - \zeta(n-\zeta)) \rho^\zeta f = O(\rho^{\zeta+1}).
$$
Since $\tkrnl{B}$ has the form $\phi_l^\zeta \times (\text{smooth})$
near the left face, and
$Q_1$ vanishes there to infinite order, we see from
(\ref{istepb}) that in fact
$$
Q_2 \in \op{}{\zeta+1,\zeta,1}
$$

It is not difficult to see that this step is continuous.
Consider the solution of (\ref{istepa}).  By linear change of coordinates, 
assume $g_p$ is the standard hyperbolic metric.  The space $T^+_pX$, on which
$N_p(\cdot)$ lives, is diffeomorphic to the interior of the fiber $S^{n+1}_{++}$. 
This fiber has boundary defining functions $\sigma_l$ and $\sigma_r$, which are
the restrictions of $\phi_l$ and $\phi_r$.  Given $Q_1
\in\op{-\infty}{}$, $N_p(Q_1)$ defines an element of $\dot\cinf(S^{n+1}_{++})$. 
We have already noted that the map
$$
\tkrnl{Q_1}\mapsto N_p(Q_1).
$$
is continuous.  Lemma 6.13 of \cite{MM} says that
$$
G_\zeta: \dot C^\infty(S^{n+1}_{++}) 
\to \sigma_l^{\zeta} \sigma_r^{\zeta} C^{\infty}(S^{n+1}_{++}),
$$
is continuous.
This gives us a function on $F_f$ which we may clearly extend into $\sxx$ so as
to preserve the continuity.  The conclusion is that the map $\tkrnl{Q_1} 
\mapsto \tkrnl{B}$ giving to the solution of (\ref{istepa}) 
may be made continuous as a map 
$$
\phgk{\infty,\infty,0} \to \phgk{\zeta,\zeta,0}.  
$$

Since $Q_2 = P_\zeta B + Q_1$, it follows from Theorem \ref{compthm} that
$\tkrnl{Q_2}$ also depends continuously on $Q_1$ and
$P_\zeta$, as an element of $\phgk{\zeta+1,\zeta,1}$.   At this stage we have
$P_\zeta(A+B) = I - Q_2$,
where the maps $g\mapsto \tkrnl{A}$, $\tkrnl{B}$ and $\tkrnl{Q}$ may be assumed
continuous in $\cnrml{-2}$, $\phgk{\zeta,\zeta,0}$, and $\phgk{\zeta+1,\zeta,1}$, 
respectively

The remaining task is to invert $I-Q_2$.
The operator $Q_2$ is compact on an appropriate weighted $L^2$ space,  
but the Neumann series for $(I-Q_2)^{-1}$ doesn't necessarily converge.
In order to make use of the Neumann series, 
we fix a particular metric $g_0$ and consider only metrics in a neighborhood of $g_0$ in
$\metr$.   By adding a term $C \in \vrop{\zeta,\zeta}$,
and adjusting the earlier construction of $A$ and $B$ accordingly, we may assume
that $R_\zeta(g_0) = (A+B+C)(g_0)$.  Of course we also assume that $C$ depends
continuously on $g$.  Thus, given seminorms on the appropriate spaces, we may
construct $A,B,C,Q$ such that
$$
P_\zeta(A+B+C) = I-Q,
$$
with all operators are continuous functions of $g$ in the appropriate spaces.
And in addition, $Q(g_0) = 0$.

\begin{lemma}\label{phgcomp}
Let $\calE$ be an index family such that 
$E_1+E_2 > n$ and $\calE\circ\calE = \calE$.  
The topology on $\op{}{\calE}$ may be defined with a family of seminorms
$\{\norm{\cdot}_a\}$, which each have the property that
\begin{equation}\label{nfcg}
\norm{A\circ B}_a \le \norm{A}_a \;\norm{B}_a,
\end{equation}
for any $A,B \in \op{}{\calE}$.
\end{lemma}

\begin{proof}
By the continuity of the composition,
$$
\op{}{\calE}\circ\op{}{\calE} \to \op{}{\calE},
$$
we can always estimate a seminorm on the right-hand side with some different
seminorms on the left-hand side.  The point we need to check is that the estimates 
on the left will require no more derivative than those on the right.  
To estimate the derivative of
$\tkrnl{A\circ B}$ by some vector field $V \in \calV_b(\sxx)$, recall how the
composition is defined in (\ref{lkcomp}).  We lift $V$ to an element of 
$\calV_b(\sxxx)$ and apply it to $\beta_{12}^* \tkrnl{A} \cdot \beta_{23}^* \tkrnl{B}$.
Clearly the result may be estimated by some combination of bounds on elements of
$\calV_b(\sxx)$ applied to $\tkrnl{A}$ and $\tkrnl{B}$, together with sup-norms
of each.  (Of course the appropriate weightings 
at the boundary must be included.)  So if $\norm{A}_a$
includes estimates the derivatives of $\tkrnl{A}$ under a spanning set of vector fields
in $\calV_b(\sxx)$ (plus estimates of the undifferentiated 
$\tkrnl{A}$), then we have (\ref{nfcg}) up to a constant.  The constant is removed by
rescaling.

By induction, we can do the same for seminorms with arbitrary numbers of derivatives.
(Our seminorms will all be norms, in fact.)

This would suffice for the conormal topology.  The polyhomogeneous conormal topology
requires in addition estimates on all coefficients in boundary expansions.  But the
argument is essentially the same.  We remove the leading terms in the boundary 
expansions (in order to study the remainders) 
by applying particular differential operators.  For example, in local coordinates
$(x,y)$, applying $x\del_x - \alpha$ to $x^\alpha f$ yields $x^\alpha \del_xf$, whose
leading term is the second coefficient in the original expansion.  
(Such operators are often used to define polyhomogeneity, as in \cite{Me92}.)
So estimates of boundary coefficients of $\tkrnl{A\circ B}$ 
may be done by lifting such operators to
$\sxxx$, and a similar argument applies.
\end{proof}

Formally, $(I-Q)^{-1} = I + Q'$, where
$$
Q' =  \sum_{j=1}^\infty Q^j.
$$
Assuming that this series converges, the resulting operator will 
be an element of $\op{}{\calE}$, where
\begin{equation*}
\begin{gathered}
E_1 = \{(\zeta+1+k,k);\;k=0,1, \dots\}, \\
E_2 = \{(\zeta+k,k);\;k=0,1, \dots\}, \\
E_3 = \{1\} \cup \{(2\zeta+k, (k^2+k-2)/2);\; k=1,2, \dots\}
\end{gathered}
\end{equation*}

Fix a seminorm $\norm{\cdot}_a$ on $\op{}{\calE}$ with the property given in
Lemma \ref{phgcomp}..  Then we define the neighborhood
$$
W_a = \{g\in \metr: \norm{Q}_a<1/2\}.
$$
The series for $Q'$ converges in $\norm{\cdot}_a$, 
uniformly for $g\in W_a$.  So in this
neighborhood we have a well-defined $Q'$ which depends continuously on
$g$ as measured by $\norm{\cdot}_a$.  We may assume that the seminorm is strong enough
to guarantee that $Q'$ is well-defined as an operator on $\rho^\delta \hfdhil$.

Now let $M = (A+B+C)(I+Q')$.  By construction, $M$ is a right inverse of 
$P_\zeta$ on $\rho^\delta \hfdhil$ for $g\in W_a$
and is a continuous function of $g$ with respect to some particular seminorm 
(related to $\norm{\cdot}_a$) on
$\op{-2}{\zeta,\zeta,0} + \op{}{(\zeta,\zeta,0)\xU \calE}$.  
We will show $M=R_\zeta$.

This will follow from the uniqueness result proved in Proposition 
\ref{uniqprop}, by the following argument. 
For $u \in \dot\cinf(X; \hfd)$, $Mu$ is certainly well-defined.
If $Q'$ actually converged in $\op{}{\calE}$, then
Theorem \ref{compthm} and Proposition \ref{opact} would imply that
$$
Mu \in \phgd{F}(X; \rho^{-(n+1)/2} \hfd),
$$ 
where the index set $F = \{ (\zeta+k,k); k = 0,1,\dots\}$.
We can't quite assume this, but by making $\norm{\cdot}_a$ strong enough,
we can at least insure that
$$
Mu  = [\rho^{\zeta} f +  O(\rho^{\zeta+\epsilon})] \;\mu_g.
$$

We also know $P_\zeta (Mu - R_\zeta u) = 0$.  The hypotheses on $\zeta$
allow us to apply Proposition
\ref{uniqprop} and conclude that $Mu = R_\zeta u$.  And since 
$\dot\cinf(X; \hfd)$ is dense in any space we would care to consider, 
this means $M = R_\zeta$. 

Since $g_0$ and $\norm{\cdot}_a$ were arbitrary, this means we can control
the continuity of $R_\zeta$ in the topology of
$\op{-2}{\zeta,\zeta,0} + \op{}{(\zeta,\zeta,0)\xU \calE}$.  
Since the topology on polyhomogeneous conormal distributions includes 
control of all asymptotic coefficients, this implies continuity in
the smaller space $\op{-2}{\zeta,\zeta,0} + \vrop{\zeta,\zeta}$
\end{proof}

\subsection{Extension to the whole plane}
From Theorem \ref{resthm} we easily deduce the following.
\begin{corollary}\label{sccor}
Assume $\re\zeta \ge \frac{n}2$, $\zeta \ne n/2$.  The map
$$
\metr \ni g \mapsto S_\zeta \in \psbop{2\zeta-n}
$$
is continuous, except at poles
\end{corollary}

\begin{proof}
The kernel of the scattering operator has the local form
$$
K_{S_\zeta}(y,y') = r^{-2\zeta} F(r,\theta,y) + G(y,y'),
$$
where $(y,y')$ are local coordinates for $\bxbx$,
$r=|y-y'|$, $\theta = (y-y')/r$, and $F$ and $G$ are smooth
in their respective variables.
These functions $F$ and $G$ are just coefficients in the boundary expansion of 
$\tkrnl{R_\zeta}$, so Theorem \ref{resthm} immediately yields the continuity of
the maps $g\mapsto F,G$ in a $\cinf$ topology.
This in turn yields the continuity of $S_\zeta$ in the topology of
$\psbop{2\zeta-n}$, except at values of $\zeta$ for which the distribution 
$r^{-2\zeta}$ has a pole.
\end{proof} 

This result may be extended to general $\zeta$ using the
relation\footnote{This idea was pointed out to me by Peter Perry.}
\begin{equation}\label{ssi} 
S_{n-\zeta} S_\zeta = I.
\end{equation}
We also need the fact that compositions
$$
\psbop{a} \circ \psbop{b} \to \psbop{a+b}
$$
are continuous, which is standard.

\begin{proof}[Proof of Theorem \ref{scthm}]
Pick a metric $g_0$ and suppose that $\re\zeta>n/2$ and $S_\zeta(g_0)$
exists and is invertible.
So there is no scattering pole at either $\zeta$ or $n-\zeta$ for $g_0$. 
That $S_\zeta(g)$ is also well-defined in a neighborhood of $g_0$ follows
from Corollary \ref{sccor}.  We will show that $S_\zeta(g)$
is also invertible for metrics in a neighborhood of $g_0$ and that this inverse is
continuous.

Let $T_\zeta(g) \in \psbop{n-2\zeta}$ be a family of pseudodifferential
parametrices for $S_\zeta(g)$ which are continuous functions of $g$
in some neighborhood of $g_0$. 
Such a family may be constructed as in the first phase of the proof of
Theorem \ref{resthm}.  
Thus we can assume
$$
S_\zeta T_\zeta = I + K_\zeta,
$$
where $g\mapsto K_\zeta \in \psbop{-\infty}$
is also continuous.  Furthermore, we can arrange that $T_\zeta(g_0) = 
S_{n-\zeta}(g_0)$, so $K_\zeta(g_0) = 0$.
From the identity (\ref{ssi}), we see that $S_{n-\zeta}$ is given by
$$
S_{n-\zeta} = T_\zeta(I+K_\zeta)^{-1},
$$
when this inverse exists.  By the same methods as in the
final phase of the proof of Theorem \ref{resthm}, $I+K_\zeta$ must be invertible
in a neighborhood of $g_0$ and will depend continuously on $g$.
\end{proof}

Finally, using Theorem \ref{scthm} we extend Theorem \ref{resthm} to the whole plane. 
\begin{proof}[Proof of Theorem \ref{mainthm}]
From Proposition \ref{rseprop} we have the relation:
\begin{equation}\label{krnz}
\krnl{R_{n-\zeta}}(q,q') = \krnl{R_\zeta}(q,q') + (n-2\zeta)\mu_g\mu'_g \int_{\bX}
E_\zeta(q, p)\; [S_{n-\zeta}E_{\zeta}](q',p) \>dh|_{\bX}.
\end{equation}
$E_\zeta$ is a polyhomogeneous conormal distribution on 
$\sxbx$, the stretched version of $X\times\bX$ defined by local polar
coordinates $r = \sqrt{x^2 + |y-y'|^2}$ and $(\eta,\theta) = (x,y-y')/r$.
It's continuity as a function of $g$ follows immediately from 
Theorem \ref{resthm} for $\re\zeta\ge n/2$.

So for $\re\zeta \ge n/2$ (and away from poles of $S_{n-\zeta}$), 
we have established the continuity of all kernels appearing 
on the right-hand side of (\ref{krnz}).  The final step is to write the 
integral over $\bX$ as a sequence of pull-backs and push-forwards, 
and apply the machinery of \cite{MeBk}.  We omit the details.
\end{proof}

\end{document}